\documentclass[11pt]{article}
\usepackage{amsthm, amssymb}
\usepackage{soul, color}
\usepackage{epsfig,rotate}
\usepackage[]{amsmath}
\usepackage[]{amsfonts}
\usepackage[]{graphicx}
\newcommand{\be}{\begin{equation}}
\newcommand{\ee}{\end{equation}}
\newcommand{\ba}{\begin{array}}
\newcommand{\ea}{\end{array}}
\newcommand{\bea}{\begin{eqnarray*}}
\newcommand{\eea}{\end{eqnarray*}}
\newcommand{\bean}{\begin{eqnarray}}
\newcommand{\eean}{\end{eqnarray}}
\newtheorem{theorem}{Theorem}
 
\newtheorem{proposition}{Proposition}
\newtheorem{lemma}{Lemma}

\def\a{{\alpha}}

\def\d{{\delta}}
\def\e{{\varepsilon}}
\def\se{{\sigma_{ess}}}
\def\lm{{\lambda^-}}
\def\lp{{\lambda^+}}

\def\S{{{\mathcal S}_D}}
\def\Ks{{{\mathcal K}^*_D}}

\def\vf{{\varphi}}

\def\R{{\mathbf R}}
\def\ds{{\displaystyle}}

\newcommand{\finproof}
{%
\mbox{}%
\nolinebreak%
\hfill%
\rule{2mm}{2mm}%
\medbreak%
\par%
}

\title{Characterization of the essential spectrum of the Neumann-Poincar\'e
operator in 2D domains with corner via Weyl sequences}
\author{
Eric Bonnetier 
\thanks{Laboratoire Jean Kuntzmann, 
Universit\'e Grenoble-Alpes, 
700 Avenue Centrale,
38401 Domaine Universitaire de Saint-Martin-d'H\`eres, France,
({\tt Eric.Bonnetier@imag.fr}).}
\and
Hai Zhang,
\thanks{Department of Mathematics, 
Hong Kong University of Science and Technology, 
Clear Water Bay, Kowloon, HK,
({\tt haizhang@ust.hk}). }
}

\begin{document}

\maketitle

\begin{abstract}
The Neumann-Poincar\'e (NP) operator naturally appears in the context of 
metamaterials as it may be used to represent the solutions of elliptic
transmission problems via potentiel theory. 
In particular, its spectral properties are closely related to the 
well-posedness of these PDE's, in the typical case where one
considers a bounded inclusion of homogeneous plasmonic metamaterial 
embedded in a homogeneous background dielectric medium.
In a recent work~\cite{PerfektPutinar_2}, 
M. Perfekt and M. Putinar have shown that the NP operator of a
2D curvilinear polygon has an essential spectrum, which depends
only on the angles of the corners. Their proof is based on quasi-conformal mappings
and techniques from complex-analysis.
In this work, we characterize the spectrum of the NP operator for a 2D domain with corners
in terms of elliptic corner singularity functions, which gives insight on the behavior of generalized eigenmodes.
\end{abstract}

%%-----------------------------------------------------------------------

\section{Introduction}
%%-----------------------------------------------------------------------

Plasmonic metamaterials are composite structures, in which some parts are made of media 
with negative indices. Their fascinating properties of subwavelength confinement 
and enhancement of electro-magnetic waves have drawned considerable interest 
from the physics and mathematics communities.
The progress in the controled production of composites with characteristic 
features of the order of optical wavelengths contributes to this 
activity, as it may enable many applications to nano-optical-mechanical systems, 
cancer therapy, neuro-science, energy and information storage and
processsing.
\medskip

From the mathematical modeling point of view, these studies have also
renewed interest in the Neumann-Poincar\'e operator, the integral operator
derived from the normal derivative of the single layer potential.
Indeed, it proves to be an interesting tool to construct, represent
and derive properties of solutions to diffusion-like equations, 
in situations where the Lax-Milgram theory does not apply, which is 
typically the case of negative index materials.
\medskip

The spectral properties of this operator have proved interesting in several 
contexts~\cite{AndoKang,AmmariCiraoloKangLeeMilton,BDT,BonnetierTriki,BonnetierTriki_2}.
They are particularly relevant to
metamaterials, as they are closely related to the existence of surface plasmons,
i.e., solutions of the governing PDE (Maxwell, Helmholtz, accoustic equations)
which are supported in the vicinity of the interfaces where the coefficients
change signs.
\medskip

To fix ideas, we consider a single inclusion $D$ made of negative
index material (typically metals, such as gold or silver at optical 
frequencies).
It is embedded in a homogeneous dielectric background medium
and we by denote $\Ks$ the associated NP operator (its precise definition 
is given in section 2). 
For particular frequencies, called plasmonic resonant frequencies,
an incident wave may excite electrons on 
the surface of the inclusion into a resonant state, that generates
highly oscillating and localized electromagnetic fields.
%The induced enhancement of the fields is quite interesting for 
%applications, as a means to control light at the subwavelength
%scale. 
For gold and silver, plasmonic resonances occur when the diameter
of the particles is small compared to the wavelength. From the modeling
point of view, one may rescale the governing Maxwell or Helmholtz equations,
with respect to particule size,
and take the limit of the resulting equations to obtain the quasi-static
regime, where only the higher-order terms of the original PDE 
remain~\cite{MayergoyzFredkinZhang,Grieser,Hai_1,Hai_2}. 
Plasmonic resonances have been investigated via layer potential
techniques in~\cite{AndoKang}--\cite{Hai_4}.
\medskip

When $D$ has a smooth boundary (say ${\cal C}^2$) the operator $\Ks$
is compact. Its spectrum is real, contained in the interval~$(-1/2,1/2]$,
and consists in a countable number of eigenvalues that accumulates to $0$.
In the context of plasmonics, domains with corners present an obvious
interest when one attemps to concentrate electro-magnetic fields,
and several authors have considered geometries where the negative
index materials are distributed in regions with 
corners~\cite{Bonnet_etal_2, Bonnet_etal_3, HelsingPerfekt, HelsingKangLim}.
When $D$ has corners, $\Ks$ is not compact~\cite{Verchota}.
In a recent work, M.-K. Perfekt and M. Putinar have shown, 
relying on the relationship betwen complex analysis and potential
theory, that the NP operator associated to a planar domain with corners 
has essential spectrum, which they characterized to be
\begin{eqnarray*}
\se(\Ks) &=& [\lambda_-, \lambda_+],
\quad \lambda_+ \;=\; -\lambda_- \;=\; \ds\frac{1}{2}(1 - \ds\frac{\a}{\pi}),
\end{eqnarray*}
where $\a$ is the most acute angle of $D$.
See~\cite{PerfektPutinar_2, PerfektPutinar_1}.
\medskip

The objective of our paper is to give an alternative derivation of the 
essential spectrum of $\Ks$ when $D$ has corners, and to establish
a close connection between the fact that $\Ks$ has essential spectrum 
and the theory of elliptic corner singularities initiated by Kondratiev
in the 1970's and developed in many directions.
See~\cite{Kondratiev} and also~\cite{Grisvard, CostabelDaugeNicaise, KozlovMazyaRossmann} 
and the many references therein.
This theory shows that the solution $u$ to an elliptic scalar equation
in a domain $O$ with corners splits as the sum $u = u_{reg}+ u_{sing}$
of a regular part $u_{reg} \in H^2(O)$ and a singular part 
$u_{sing} \in H^1(O) \setminus H^2(O)$, locally around each corner. 
Up to a scaling factor, the expression of the latter part,
which we call `singularity function', only depends 
on the geometry of the corner, and on the nature of the boundary conditions.
In the case of a transmission problem, it depends on the angle and on
the contrast in material coefficients. 
Typically, $u_{sing}$ is a non-trivial solution of a
homogeneous problem for the associated operator in the infinite domain
obtained by zooming around the vertex of the corner.
For a transmission problem in 2D, it has the form
\begin{eqnarray} \label{form_using}
u_{sing} &=& C r^\eta \vf(\theta),
\end{eqnarray}
where $(r,\theta)$ denote the polar coordinates with orgin at the
vertex of the corner under consideration. The exponent $\eta$
is the root of a dispersion relation, and $\vf$ is a smooth function
(or piecewise smooth in the case of a transmission problem).
\medskip

This paper is organized in the following way.
Section~2 of the paper describes the setting and notations and reviews
useful facts about the NP operator. In Section~3, we study how elliptic
corner singularity functions depend on the conductivity contrast.
In the very interesting papers~\cite{Bonnet_etal_1, Bonnet_etal_2, Bonnet_etal_3}, 
it is shown that
functions of the form~(\ref{form_using}) only exist when the conductivity
contrast $\lambda$ lies outside a critical interval $[\lambda_-, \lambda_+]$.
When $ \lambda \in [\lambda_-, \lambda_+]$, the elliptic corner singularity 
functions still have the form $u_{sing}$ but their expression involves a
complex exponent $\eta$. 
In~\cite{Bonnet_etal_3}, the use of the Mellin transform converts the search of
these singular functions to that of propagative mode in an infinite
wave-guide. These functions are called plasmonic black-hole waves, 
reflecting the fact that they are not in the
energy space $H^1(\Omega)$. 
In Section~4, we show that the critical interval is contained in the essential 
spectrum $\sigma_{ess}(\Ks)$, by generating singular Weyl sequences~\cite{BirmanSolomjak} 
using the singularity functions.
In Section~5, the reverse inclusion is proved. In particular, we use a construction
inspired by~\cite{HoaiMinh} to transform, around the vertex of the corner, 
the PDE with sign changing conditions into a system of PDE's defined in the
inhomogeneity only, that satisfies complementing boundary conditions in the
sense of Agmon, Douglis and Nirenberg, and for which we prove well-posedness.

%%-----------------------------------------------------------------------

\section{The Neumann-Poincar\'e operator and the Poin\-car\'e variational operator}
%%-----------------------------------------------------------------------

Throughout the text, $\Omega \subset \R^2$ denotes a bounded 
open set with smooth boundary, that strictly contains a connected inclusion $D$.
We assume that $\partial D$ is smooth, except for one corner point, 
of angle $\a, 0 < \a <Pi$, located at the origin.
We assume that for some $R_0 > 0$, 
\begin{eqnarray} \label{A1}
D \cap B_{R_0} &=& \{ x = (r\cos(\theta), r \sin(\theta)),
0 \leq r < R_0, |\theta| < \a/2 \},
\end{eqnarray}
where, for any $\rho > 0$,
$B_{\rho}$ denotes the ball of radius $\rho$ centered at $0$.
The space $H^1_0(\Omega)$ is equipped with the following inner product and
associated norm
\begin{eqnarray*}
<u,v>_{H^1_0} \;=\; \ds\int_\Omega \nabla u \cdot \nabla v \,dx,
&\quad&
||u||_{H^1_0} \;=\; 
\left(\ds\int_\Omega |\nabla u|^2 \,dx
\right)^{1/2}.
\end{eqnarray*}
Our work concerns the following diffusion equation: given a function 
$f~\in~L^2(\Omega)$, we seek $u$ such that
\begin{equation} \label{eq_cond}
\left\{ \begin{array}{ccll}
-\textrm{div}(a(x) \nabla u(x)) &=& f &\textrm{in}\; \Omega,
\\
u(x) &=& 0 & \textrm{on}\; \partial \Omega,
\end{array} \right.
\end{equation}
where the conductivity $a$ is piecewise constant
\begin{eqnarray} \label{def_a}
a(x) &=&
\left\{ \begin{array}{ccll}
k \in {\mathbf C}& \quad x \in D,
\\
1 & \quad x \in \Omega \setminus \overline{D}.
\end{array} \right.
\end{eqnarray}
\medskip

It is well known that when $k$ is strictly positive, or when $k \in \mathbf{C}$
and $Im(k)~\neq~0$, this problem has a unique solution in $H^1(\Omega)$, and that
\begin{eqnarray*}
||u||_{H^1_0} &\leq& C(k) \, ||f||_{L^2},
\end{eqnarray*}
for some constant $C(k) > 0$ that depends on $k$.
\medskip

Let $P(x,y)$ denote the Poisson kernel associated to $\Omega$, defined by
\begin{eqnarray*}
P(x,y) &=& G(x,y) + R_x(y), \quad x,y \in \Omega,
\end{eqnarray*}
where $G(x,y)$ denotes the free space Green function
\begin{eqnarray*}
G(x,y) &=& \ds\frac{1}{2\pi} \ln|x-y|,
\end{eqnarray*}
and where $R_x(y)$ is the smooth solution to
\[
\left\{ \begin{array}{clcl} 
\Delta_y R_x(y) &=& 0 & y \in \Omega,
\\
R_x(y) &=& -G(x,y) & y \in \partial \Omega.
\end{array} \right.
\]
With the Poisson kernel, we define the single layer potentials
$\S\vf \in L^2(\partial D)$ of a function $\vf \in L^2(\partial D)$ by
\begin{eqnarray*}
\S\vf(x) &=&
\int_{\partial D} P(x,y) \vf(y)\,ds(y),
\quad x \in D \cup (\Omega \setminus \overline{D}).
\end{eqnarray*}
It is well known~\cite{Folland,Verchota} that $\S\vf$ is harmonic in $D$ and in 
$\Omega \setminus \overline{D}$, continous in $\overline{\Omega}$, and that 
its normal derivatives satisfy the Plemelj jump conditions
\begin{eqnarray} \label{jump_S}
\ds\frac{ \partial \S \vf}{\partial \nu}|^\pm(x)
&=&
(\pm\ds\frac{1}{2}I + \Ks)\vf(x),\quad x \in \partial D.
\end{eqnarray}
where $\Ks$ is the Neumann-Poincar\'e operator, defined by
\begin{eqnarray*}
\Ks\vf(x) &=&
\ds\int_{\partial D} \ds\frac{\partial P}{\partial \nu_y}(x,y)
\vf(y) \,ds(y).
\end{eqnarray*}
It is shown in~\cite{Coifman_etal} that this definition makes sense
for Lipschitz domains, and in that case, the operator $\Ks$ is
continuous from $L^2(\partial D) \rightarrow L^2(\partial D)$,
which extends as an operator 
$H^{-1/2}(\partial D) \rightarrow H^{1/2}(\partial D)$.
\medskip

The solution $u$ to~(\ref{eq_cond}) can then be represented in the form
\begin{eqnarray} \label{repr_u}
u(x) &=&
\S\vf(x) + H(x),
\end{eqnarray}
where the harmonic part is given by
\begin{eqnarray*}
H(x) &=& \int_{\Omega} P(x,y) f(y)\, ds(y).
\end{eqnarray*}
The jump conditions~(\ref{jump_S}), constrain the layer potential 
$\vf \in H^{-1/2}(\partial D)$ to satisfy the integral equation
\begin{eqnarray*} \label{eq_int}
(\lambda I - \Ks) \vf(x) &=&
\partial_\nu H_{|\partial D}(x),
\quad x \in \partial D.
\end{eqnarray*}
\medskip

We also introduce the Poincar\'e variational operator 
$T_D~: H^1_0(\Omega) \rightarrow H^1_0(\Omega)$, 
defined for $u \in H^1_0(\Omega)$ by
\begin{eqnarray} \label{def_T}
\forall \; v \in H^1_0(\Omega), \quad
\ds\int_\Omega \nabla T_D u \cdot \nabla v \,dx
&=& 
\ds\int_D \nabla u \cdot \nabla v \,dx.
\end{eqnarray}
Some of its properties are described in the following proposition
(see~\cite{BDT} for a proof).

%%-----------------------------------------------------------------------
\begin{proposition} \label{prop_TD}
The operator $T_D$ is bounded, selfadjoint, and satisfies $||T_D||~=~1$.
Moreover, 
\begin{itemize}
\item[(i)]
Its spectrum $\sigma(T_D)$ is contained in the interval $[0,1]$.

\item[(ii)] Its kernel, the eigenspace associated to $\beta = 0$, is
\begin{eqnarray*}
Ker(T_D) &=&
\{ u \in H^1_0(\Omega), u = const \;\textrm{on}\; D \}.
\end{eqnarray*}

\item[(iii)] $1 \in \sigma(T_D)$ and the associated eigenspace is
\begin{eqnarray*}
Ker(I-T_D) &=&
\{ u \in H^1_0(\Omega), u = 0 \;\textrm{in}\; \Omega \setminus \overline{D} \},
\end{eqnarray*}
(and thus, can be identified with $H^1_0(D)$).

\item[(iv)]
The space $H^1_0(\Omega)$ decomposes as
\begin{eqnarray*}H^1_0(\Omega) &=&
Ker(T_D) \oplus Ker(I-T_D) \oplus {\mathcal H},
\end{eqnarray*}
where ${\mathcal H}$ is the closed subspace defined by
\begin{eqnarray*}
{\mathcal H} &=&
\{ u \in H^1_0(\Omega), \Delta u = 0 \;
\textrm{in}\; D \cup (\Omega \setminus \overline{D}),
\ds\int_{\partial D} \ds\frac{\partial u^+}{\partial \nu} \, ds = 0 \}.
\end{eqnarray*}
\end{itemize}

\end{proposition}
%%-----------------------------------------------------------------------

The space $\mathcal{H}_s = \mathcal{H} \oplus Ker(T_D)$ is the space
of single layer potentials. It is isomorphic to
\begin{eqnarray*}
H^{-1/2}_0(\partial D)
&=& \{ \vf \in H^{1/2}(\partial D), \;
<\vf,1>_{H^{-1/2},H^{1/2}} = 0 \}.
\end{eqnarray*}
This latter space is equipped with the inner product
\begin{eqnarray} \label{def_scalp}
<\vf,\psi>_S &=& -\ds\int_{\partial D} \vf S \psi\,d\sigma
\end{eqnarray}
for which the operator 
$\Ks~: H^{-1/2}_0(\partial D) \rightarrow H^{-1/2}_0(\partial D)$
is self-ajdjoint as a result of the Calder\'on identity~\cite{KhavisonPutinarShapiro}.
We denote by $||\cdot||_S$ the associated norm.
In particular, if $u, v \in \mathcal{H}_s$ are such that 
$u = S_D \vf, v = S_D \psi$, then the jump conditions~(\ref{jump_S})
and integration by parts show that
\begin{eqnarray} \label{scalp_uv}
\ds\int_{\Omega} \nabla u \cdot \nabla v
&=& <\vf,\psi>_S.
\end{eqnarray}
\medskip

When the domain $D$ has a ${\cal C}^2$ boundary, the Poincar\'e-Neumann 
operator $\Ks~: H^{-1/2}_0(\partial D) \rightarrow H^{-1/2}_0(\partial D)$ 
is compact. Its spectrum $\sigma(\Ks)$ is contained in $[-1/2,1/2]$, 
and consists of a sequence of real eigenvalues that accumulates to 0.
In this case, $\sigma(\Ks)$ is directly related to $\sigma(T_D)$.
Indeed, if $u \in H^1_0(\Omega)$ and $\beta \in \R, \beta \neq 1$, satisfy
$T_D u = \beta u$, it follows from~(\ref{def_T}) that
\begin{eqnarray*}
\forall \; v \in H^1_0(\Omega), \quad
\beta \ds\int_{\Omega \setminus D} \nabla u \cdot \nabla v \,dx
+ (\beta -1) \ds\int_D \nabla u \cdot \nabla v \,dx
&=& 0,
\end{eqnarray*}
so that $u$ is a non-zero solution to
\begin{equation} \label{pde_0rhs}
\left\{ \begin{array}{ccll}
\textrm{div}(a(x)\nabla u(x)) &=& 0 & \textrm{in}\; \Omega,
\\
u(x) &=& 0 & \textrm{on}\; \partial \Omega,
\end{array} \right.
\end{equation}
where the conductivity $a$ equals $\beta$ in $\Omega \setminus \overline{D}$
and $(\beta - 1)$ in $D$. Expressing $u$ in the form $u = \S\vf$ yields
yields the integral equation
\begin{eqnarray*}
(\lambda I - \Ks)\vf(x) &=& 0, \quad x \in \partial D,
\end{eqnarray*}
where $\lambda = 1/2 - \beta$ is thus an eigenvalue of $\Ks$.
It follows that
\begin{eqnarray*}
\sigma(T_D) &=& (1/2 - \sigma(\Ks)) \cup \{0,1\}
\end{eqnarray*}
\medskip

As recalled above, when $D$ is a domain with corners, 
$\sigma(\Ks)$ contains an interval of essential spectrum~\cite{PerfektPutinar_2}.
We have
%%-----------------------------------------------------------------------
\begin{proposition} \label{prop_spect}
The essential spectra of $T_D$ and $\Ks$ are related by
$\sigma_{ess}(T_D) = 1/2 - \sigma_{ess}(\Ks)$.
\end{proposition}
%%-----------------------------------------------------------------------
{\bf Proof:}
Let $\lambda \in \sigma_{ess}(\Ks)$. By definition, there exists
a singular Weyl sequence, i.e.,
a sequence of functions $(\vf_\e) \subset H^{-1/2}_0$ such that
\[
\left\{ \begin{array}{ccll}
(\lambda I - \Ks)\vf_\e &\rightarrow& 0 &\textrm{strongly in}\; H^{-1/2}_0,
\\
||\vf_\e||_S &=& 1,
\\
\vf_\e & \rightharpoonup & 0 &\textrm{weakly in}\; H^{-1/2}_0.
\end{array} \right.
\]
Let $\beta = 1/2 - \lambda$ and $u_\e = S_D \vf_\e \in \mathcal{H}_S$. 
%%and it follows from~(\ref{jump_S}) that
%%\begin{eqnarray} \label{rel_ue}
%%\vf_\e \;=\; \partial_\nu u_\e \lvert^+ - \partial_\nu u_\e \lvert^-,
%%&\textrm{and}&
%%\Ks \vf \;=\;  \partial_\nu u_\e \lvert^+ + \partial_\nu u_\e \lvert^-.
%%\end{eqnarray}
Let $v \in \mathcal{H}_S$ so that $v = S_D \psi$ for some 
$\psi \in H^{-1/2}_0(\partial D)$. It follows from~(\ref{scalp_uv}) that
\begin{eqnarray*}
\ds\int_\Omega \nabla u_\e \cdot \nabla v
&=& 
<\vf_\e, \psi> \;\rightarrow\; 0.
\end{eqnarray*}
This equality also holds for $v \in Ker(I - T_D)$ since this subspace
is orthogonal to $\mathcal{H}_S$, and thus
\begin{eqnarray} \label{cond_1}
u_\e &\rightharpoonup& 0 \quad \textrm{weakly in}\; H^1_0(\Omega).
\end{eqnarray}
\medskip

Additionally, invoking~(\ref{scalp_uv}) again, we see that
\begin{eqnarray} \label{cond_2}
\ds\int_\Omega |\nabla u_\e|^2 &=& <\vf_\e,\vf_\e>_S \;=\; 1.
\end{eqnarray}
\medskip

Finally, we compute, for $v = S_D\psi \in \mathcal{H}_S$,
\begin{eqnarray*}
\ds\int_\Omega
\nabla \left( (\beta I - T_D)u_\e \right) \cdot \nabla v
&=&
\ds\int_\Omega \beta \nabla u_\e  \cdot \nabla v
\;-\;
\ds\int_D \nabla u_\e  \cdot \nabla v
\\
&=&
\ds\int_{\Omega \setminus D} \beta \nabla u_\e  \cdot \nabla v
\;+\;
\ds\int_{D} (\beta - 1) \nabla u_\e  \cdot \nabla v
\\
&=&
- \beta \ds\int_{\partial D} \partial_\nu u_\e\lvert^+ v
\;+\;
(\beta-1) \ds\int_{\partial D} \partial_\nu u_\e\lvert^- v.
\end{eqnarray*}
Inserting~(\ref{jump_S}) in place of the normal derivatives of $u_\e$
we see that 
\begin{eqnarray*}
\ds\int_\Omega
\nabla \left( (\beta I - T_D)u_\e \right) \cdot \nabla v
&=& <(\lambda I - \Ks)\vf_\e,\psi>_S
\\
&\leq& ||(\lambda I - \Ks)\vf_\e||_S \, ||\psi||_S.
\end{eqnarray*}
It follows that 
\begin{eqnarray} \label{cond_3}
||(\beta I - T_D)u_\e||_{H_1} &\leq& ||(\lambda I - \Ks)\vf_\e||_S
\;\rightarrow 0.\end{eqnarray}
we conclude from~(\ref{cond_1}--\ref{cond_3}) that $u_\e$ is a singular
Weyl sequence associated to $\beta$, so that $\beta \in \sigma_{ess}(T_D)$.
The same argument proves the reverse inclusion 
$\sigma_{ess}(T_D) \subset \left( 1/2 - \sigma_{ess}(\Ks) \right)$.
\finproof

\section{Corner singularity functions}
%%-----------------------------------------------------------------------

Elliptic corner singularities have been the subject of much research since
the pionneering works of Kondratiev~\cite{Kondratiev}, Grisvard~\cite{Grisvard}
(see also~\cite{CostabelDaugeNicaise, KozlovMazyaRossmann}).
Essentially, the theory focuses on the regularity of solutions
to elliptic PDEs near a corner of the domain,
or in the case of a transmission problem such as~({\ref{eq_cond}),
near a corner of the interface between several phases. 
The following is a typical statement:

%%-----------------------------------------------------------------------
\begin{theorem} \label{thm_corner}
Let $k > 0$. The solution $u \in H^1_0(\Omega)$ to~(\ref{eq_cond}) decomposes as
\begin{eqnarray*}
u &=& u_{sing} + u_{reg},
\end{eqnarray*}
where $u_{reg} \in H^2(\Omega)$ and where $u_{sing}$ has the form
\begin{eqnarray} \label{def_using}
u_{sing}(x) &=& r^\eta \vf(\theta) \zeta(x), \quad x \in \Omega.
\end{eqnarray}
Here $x = (r \cos(\theta), r \sin(\theta))$ in polar coordinates,
$\zeta$ is a smooth cut-off function,
such that, for some $s > 0$
\begin{eqnarray*}
\zeta(x) &=& \left\{ \begin{array}{ll}
1 & |x| < s,
\\
0 & |x| > 2s.
\end{array} \right.
\end{eqnarray*}
Moreover, for some constant $C = C(\a, k)$, the following estimate holds
\begin{eqnarray} \label{est_corner}
||u_{sing}||_{H^1(\Omega)}
+ ||u_{reg}||_{H^2(\Omega)}
&\leq&
C \, \left( ||u||_{H^1(\Omega)} + ||f||_{L^2(\Omega)} \right).
\end{eqnarray}
\end{theorem}
%%-----------------------------------------------------------------------

These results have been derived only in the case of elliptic 
media, i.e., when $k > 0$. In the context of plasmonic metamaterials, it is
natural to try to extend them to complex values of $k$.
To our best knowledge,
the first steps in this direction
have been obtained in~\cite{Bonnet_etal_1,Bonnet_etal_2} and
concern the existence of singularity functions of the form~(\ref{def_using}).

%%-----------------------------------------------------------------------

\subsection{Regular corner singularity functions}

In this paragraph, we investigate whether one can define singular functions 
such as~(\ref{def_using}) when $k$ may also take negative values.
More precisely, we seek $H^1_{loc}(\R^2)$ solutions to
\begin{eqnarray} \label{eq_R2}
\textrm{div}(a(x) \nabla u(x)) &=& 0 \quad \textrm{in}\; \R^2,
\end{eqnarray}
of the form 
\begin{eqnarray} \label{form_using}
u(x) = r^\eta \vf(\theta), \eta \in \R,
\end{eqnarray}
when the conductivity $a(x)$ is defined in the whole of $\R^2$ by
\begin{eqnarray} \label{def_a2}
a(x) &=& \left\{ \begin{array}{ll}
k & |\theta| < \a/2,
\\
1 & \textrm{otherwise}.
\end{array} \right.
\end{eqnarray}

Since we are only interested in singular solutions which belong to 
$H^1(\Omega) \setminus H^2(\Omega)$, we may restrict $\eta$
to lie in $(0,1)$.
As $u$ is harmonic in each sector $|\theta| < \a/2$ and
$\a/2 < \theta < 2 \pi - \a/2$, it follows that $\vf$ has the form
\begin{eqnarray} \label{expr_u}
\vf = \left\{ \begin{array}{ll}
a_1 \cos(\eta (\theta+ \a/2)) + b_1 \sin(\eta (\theta+\a/2)) 
& \textrm{if}\; -\a/2 < \theta < \a/2,
\\
a_2 \cos(\eta (\theta+\a/2)) + b_2 \sin(\eta (\theta+\a/2)) 
& \textrm{if}\; \a/2 < \theta < 2\pi - \a/2
\end{array} \right.
\end{eqnarray}
for some $a_i, b_i, i=1,2$. Expressing the continuity of $u$ and of
$a(x) \partial_\nu u$ across the interfaces, shows that a non-trivial solution
exits if and only if the following dispersion relation is satisfied
\begin{eqnarray*}
\textrm{det}\left( \begin{array}{cccc}
1 & 0 & - \cos(2 \pi \eta) & - \sin(2\pi \eta)
\\
\cos(\a \eta) & \sin(\a \eta) & - \cos(\a \eta) & - \sin(\a \eta)
\\
0 & k & \sin(2 \pi \eta) & - \cos(2\pi \eta)
\\
-k \sin(\a \eta) & k \cos(\a \eta) & - \sin(\a \eta) & - \cos(\a \eta)
\end{array} \right)
&=& 0,
\end{eqnarray*}
which, after elementary manipulations, can be rewritten in the form
\begin{eqnarray} \label{disp_rel}
\ds\frac{2k}{k^2 + 1}
&=&
\ds\frac{\sin(\a \eta) \sin((2\pi - \a)\eta)}
{1 - \cos(\a \eta) \cos((2\pi - \a)\eta)}
\;=:\; F(\eta,\a).
\end{eqnarray}
\medskip

A Taylor expansion around the values $\eta = 0$ shows that $F(\cdot,\a)$ 
can be extended by continuity to a function defined on the whole of $[0,1]$
by setting
\begin{eqnarray*}
F(0,\a) &=& \ds\frac{-2 \a (2\pi-\a)}{\a^2 + (2\pi-\a)^2}.
\end{eqnarray*}
By solving 
\[
\ds\frac{2k}{k^2 + 1}= F(0,\a),
\]
we obtain two solutions
\begin{eqnarray} \label{def_kpm} 
k_+ \;=\; \ds\frac{-(2\pi - \a)}{\a}, &&
k_- \;=\; \ds\frac{-\a}{2\pi-\a}.
\end{eqnarray}
Additionnally, it is easy to check that $|F(\eta,\a)| \leq 1$ and 
\begin{eqnarray*}
\partial_\eta F
&=&
\ds\frac{\cos((2\pi-\a)\eta)-\cos(\a \eta))
\left[ a\sin((2\pi-\a)\eta) - (2\pi-\a)\sin(\a \eta) \right] }
{\left[ 1 - \cos(\a \eta) \cos((2\pi - \a)\eta) \right]^2}.
\end{eqnarray*}
We show below that $F(\cdot,\a)$ is thus strictly increasing,
and note that $\partial_\eta F(0,\a)=\partial_\eta F(1,\a)=0$.

%%----------------------------------------------------------------
\begin{lemma} \label{prop_dF}
For any $0 < \a < \pi$ and $0 \leq \eta \leq 1$, the following inequalities hold
\begin{eqnarray}
\cos((2\pi-\a)\eta)-\cos(\a \eta) &<& 0,
\label{ineq_1}
\\
a\sin((2\pi-\a)\eta) - (2\pi-\a)\sin(\a \eta) &<& 0.
\nonumber
\end{eqnarray}
\end{lemma}
%%----------------------------------------------------------------
\medskip

{\bf Proof:}
To prove the first inequality, we first note that $\a < (2\pi - \a)$
so that $\a\eta < (2\pi - \a)\eta$. If $(2\pi - \a)\eta \leq \pi$,
then~(\ref{ineq_1}) follows from the monotonicity of the cosine function on $[0,\pi]$.
If $(2\pi - \a)\eta > \pi$, then 
\begin{eqnarray*}
\cos((2\pi-\a)\eta) &=& \cos(\pi - \beta),
\quad\textrm{with}\;(2\pi -\a)\eta =: \pi + \beta. 
\end{eqnarray*}
Noticing that
\begin{eqnarray*}
\a\eta &\leq& \a \eta + 2\pi(1 - \eta)
\;=\; \pi - \beta \;<\; \pi,
\end{eqnarray*}
we infer that $\cos(\pi - \beta) < \cos(\a\eta)$, which yields the result.
\medskip

The second inequality follows from the fact that
\begin{eqnarray*}
\lefteqn{
\partial_\eta \left[ a\sin((2\pi-\a)\eta) - (2\pi-\a)\sin(\a \eta) \right]}
\\
&=&
\a(2\pi - \a) \left[ \cos((2\pi- \a)\eta) - \cos(\a \eta) \right],
\end{eqnarray*}
which according to~(\ref{ineq_1}) is negative.
\finproof

As a consequence of~(\ref{disp_rel}), we obtain 

\begin{proposition}
Singular
solutions in $H^1_{loc}(\R^2)$ of the form~(\ref{form_using}) exists for the equation 
(\ref{eq_R2})only when 
$k \in (-\infty,k_+) \cup (k_-, +\infty)$, see Figure~\ref{fig_1}.
In terms of the contrast $\lambda = \frac{k+1}{2(k-1)}$
this condition is equivalent to
\begin{eqnarray*}
\lambda &\notin& [\lambda_-,\lambda_+]
\;:=\; [-\ds\frac{1}{2}(1 - \ds\frac{\a}{\pi}),\ds\frac{1}{2}(1 - \ds\frac{\a}{\pi})].
\end{eqnarray*}
In other words, singular solutions of the form~(\ref{form_using}) only exist when
$\lambda = \ds\frac{k+1}{2(k-1)}$ is not in $\sigma_{ess}(\Ks)$.
\end{proposition}

%%----------------------------------------------
\begin{figure}[htbp]\label{fig_1}
\begin{center}
\includegraphics[angle=0,width=50mm]{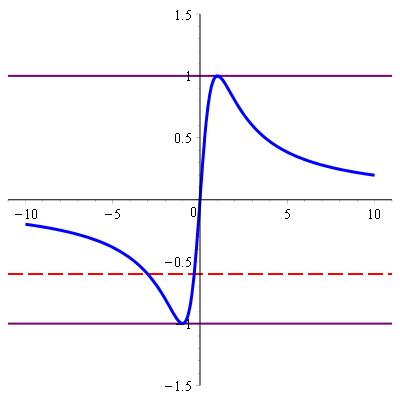}
\includegraphics[angle=0,width=50mm]{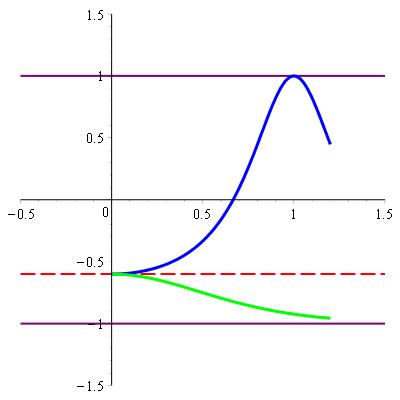}
\caption{Left: Plot of the function $k \rightarrow 2k/(k^2+1)$.
Right: Plot of $\eta \rightarrow F(\eta,\a)$ (blue), and of
$\xi \rightarrow \tilde{F}(\xi,\a)$ (green), for $\a = \pi/2$. The dotted line
indicates the value of $2k/(k^2+1)$ below which the dispersion relation 
has no solution $\eta \in \R$.}
\end{center}
\end{figure}

%%----------------------------------------------
\subsection{Singular corner singularity functions} \label{sec_3.2}

We now construct local singular solutions when $k \in [k_+,k_-]$.
By this we mean functions which satisfy the PDE~(\ref{eq_R2}),
but which may only be in $H^1_{loc}(\R^2 \setminus \{0\})$.
To this end, we seek $u(x) = r^{\eta} \vf(\theta)$, with $\vf$
in the form~(\ref{expr_u}), but assume now
that $\eta \in {\mathbf C}$. The same algebra leads to the same dispersion
relation~(\ref{disp_rel}). In particular if we restrict $\eta$
to be a pure imaginary number, $\eta = i \xi$, this relation takes
the form
\begin{eqnarray*}
\ds\frac{2k}{k^2 + 1}
&=&
\ds\frac{\sinh(\a \xi) \sinh((2\pi-\a)\xi)}
{1 - \cosh(\a \xi) \cosh((2\pi- \a)\xi)}
\;=:\; \tilde{F}(\xi,\a).
\end{eqnarray*}
It is easy to check that the function $\xi \rightarrow \tilde{F}(\xi,\a)$
can be extended by continuity at $\xi =0$ by setting
\begin{eqnarray*}
\tilde{F}(0,\a) &=& \ds\frac{2 \a (2\pi-\a)}{\a^2 + (2\pi-\a)^2}
\;=\; \ds\frac{2 k_\pm}{k_\pm^2 + 1}
\end{eqnarray*}
(and we note that $\tilde{F}(0,\a) = F(0,\a)$).
In addition, we compute
\begin{eqnarray*}
\partial_\xi \tilde{F}
&=&
\ds\frac{(\cosh((2\pi-\a)\xi) - \cosh(\a\xi))
\left[(2\pi - \a)\sinh(\a\xi) - \a\sinh((2\pi-\a)\xi) \right]}
{\left[ 1 - \cosh(\a \xi) \cosh((2\pi- \a)\xi) \right]^2}.
\end{eqnarray*}
Just as in Proposition~\ref{prop_dF}, one can show that
for any $\xi > 0$ and $0 < \a < \pi$, the product of the two factors in the above
numerator is negative, so that $\tilde{F}(\cdot,\a)$ is strictly decreasing
on $\R^+$ and its range is equal to $[F(0,\a),-1)$, see Figure~\ref{fig_1}.
We also note that $\lim_{\eta \to \infty}\tilde{F}(\eta,\a) = -1$,
so that the value $-1$ (which corresponds to $k=-1$) is never attained.
Summarizing, we have shown that

%%-----------------------------------------------------------
\begin{proposition} \label{prop_BH}
For any value of $\lambda \in (\lambda_-,\lambda_+), \lambda \neq 0$, there exists
$\xi > 0$ and a function $u(x) = r^{i\xi}\vf(\theta)$, which is
a local solution to $\textrm{div}(a(x)\nabla u(x)) = 0$,
where $a$ is defined by~(\ref{def_a2}), with
$\lambda = \frac{k+1}{2(k-1)}$.
\end{proposition}
%%-----------------------------------------------------------

%%-----------------------------------------------------------------------

\section{Construction of singular Weyl sequences}

In this section we prove
%%---------------------------------------------------
\begin{theorem} \label{thm_singWseq}
The set $[\lambda_-, \lambda_+]$ is contained in $\sigma_{ess}(T_D)$.
\end{theorem}
%%---------------------------------------------------

{\bf Proof:}

Since $\sigma_{ess}(\Ks)$ is a closed set, it is sufficient to show
that $(\lambda_-, \lambda_+) \setminus \{0\} \subset \sigma_{ess}(T_D)$.
We proceed as follows: We consider 
$\lambda \in (\lambda_-, \lambda_+), \lambda \neq 0$, 
and show that $\beta = 1/2-\lambda \in \sigma_{ess}(T_D)$
by constructing a singular Weyl sequence, i.e.,
a sequence of functions $u_\e \in H^1_0(\Omega)$, such that
\begin{equation} \label{def_singW}
\left\{ \begin{array}{ccll}
||u_\e||_{H^1_0} &=& 1,
\\
(\beta I - T_D)u_\e &\to& 0 & \textrm{strongly in}\; H^1_0(\Omega),
\\
u_\e &\rightharpoonup& 0 & \textrm{weakly in}\; H^1_0(\Omega).
\end{array} \right.
\end{equation}
\medskip
 
According to Proposition~\ref{prop_BH}, there exists $\xi > 0$ and 
coefficients $a_1, b_2, a_2, b_2 \in \mathbf{C}$, not all equal to $0$,
such that the function 
\begin{eqnarray} \label{def_usingsing}
u(x) = Re(r^{i\xi}) \phi(\theta)=
\left\{ \begin{array}{l}
Re(r^{i\xi}) \left[ a_1 \cos(i\xi(\theta+\a/2)) + b_1 \sin(i\xi(\theta+\a/2))\right] 
\\
\quad\quad \textrm{if}\;
-\a/2 < \theta < \a/2,
\\
\\
Re(r^{i\xi}) \left[ a_2 \cos(i\xi(\theta+\a/2)) + b_2 \sin(i\xi(\theta+\a/2))\right] 
\\
\quad\quad \textrm{otherwise},
\end{array} \right.
\end{eqnarray}
is harmonic in $(D \cap B_{R_0}) \setminus \{0\}$ and 
in $\left((\Omega \setminus \overline{D}) \cap B_{R_0}\right) \setminus \{0\}$,
and satisfies the transmission conditions at 
the interfaces $\theta = \pm \a/2$.
\medskip

Let $r_0 < R_0/2$ and let $\chi_1, \chi_2 : \R^+ \rightarrow [0,1]$ 
denote two smooth cut-off functions, such that for some constant $C > 0$
\[ 
\left\{ \begin{array}{lcllclcll}
\chi_1(s) &=& 0 & |s| \leq 1, &\hspace*{10mm}& \chi_2(s) &=& 0 & |s| \geq 2r_0,
\\
\chi_1(s) &=& 1 & |s| \geq 2, &\hspace*{10mm}& \chi_2(s) &=& 1 & |s| \leq r_0,
\\
|\chi_1^\prime(s)| &\leq& C, &&& |\chi_2^\prime(s)| &\leq& C.
\end{array} \right.
\]
We set $\chi_1^\e(r) = \chi_1(r/\e)$, and define 
\begin{eqnarray} \label{def_ue}
u_\e(x) &=& s_\e \chi^\e_1(r) \chi_2(r) u(x), \quad x \in \Omega.
\end{eqnarray}
%%where $u$ as in~(\ref{def_usingsing}).
\medskip

The function $u$ is not in $H^1$ as its gradient blows up
like $r^{-1}$ near the corner, consequently
\begin{eqnarray*}
m_\e &:=& \ds\int_\e^{r_0}\int_0^{2\pi} |\nabla u(x)|^2 \,r dr d\theta
\;\rightarrow\; \infty
\quad \textrm{as}\; \e \to 0.
\end{eqnarray*}
We choose $s_\e$ in~(\ref{def_ue}) so that $||u_\e||_{H^1_0} = 1$, in other words
\begin{eqnarray*}
s_\e^{-2} &=&
\ds\int_\e^{2\e}\int_0^{2\pi} |u \nabla \chi_1^\e + \chi_1^\e \nabla u|^2
\;+\;
m_\e 
\;+\;
\ds\int_{r_0}^{2r_0}\int_0^{2\pi} |u \nabla \chi_2 + \chi_2 \nabla u|^2
\\
&=:& J_1 + m_\e + J_2.
\end{eqnarray*}
The term $J_2$ is independent of $\e$ and is $O(1)$, and in particular 
\begin{eqnarray*}
J_2 & = & o(m_\e) \quad \textrm{as}\; \e \to 0.
\end{eqnarray*}
The other term can be estimated as follows
\begin{eqnarray} \label{est_J1}
J_1 &=&
\ds\int_\e^{2\e}\int_0^{2\pi} \left|
\ds\frac{r^{i\xi}+r^{-i\xi}}{2} \vf(\theta) \chi_1^\prime(r/\e)/\e + 
i\xi \ds\frac{r^{i\xi -1}-r^{-i\xi-1}}{2} \vf(\theta) \chi_1(r/\e) \right|^2
\nonumber \\
&&\hspace*{16mm}
+ \left|\ds\frac{r^{i\xi -1}-r^{-i\xi-1}}{2} \vf^\prime (\theta) \chi_1(r/\e) \right|^2
\,r dr d\theta
\\
&\leq&
C \, \ds\int_0^{2\pi} \left( |\vf(\theta)|^2 + |\vf^\prime(\theta)|^2 \right)d\theta
\,
\ds\int_\e^{2\e} \left( ||\chi_1^\prime||^2_{\infty}/\e^2 + r^{-2}||\chi_1||_\infty \right)
\,rdr
\nonumber \\
&\leq&
C\, \ds\int_0^{2\pi} \left( |\vf(\theta)|^2 + |\vf^\prime(\theta)|^2 \right)d\theta
\,
\left(
3/2 ||\chi_1^\prime||_\infty^2 + (\ln(2\e) - \ln(\e)) ||\chi_1||_\infty
\right).
\nonumber 
\end{eqnarray}
Since $\vf$ is independent of $\e$, we see that
\begin{eqnarray*}
J_1 &=& O(1) \;=\; o(m_\e), \quad \textrm{as}\; \e \to 0,
\end{eqnarray*}
and so $s_\e \sim m_\e^{-1/2} \to 0$.
\medskip

We next show that $||(\beta I - T_D)u_\e||_{H^1} \to 0$. Indeed,
let $v \in H^1_0(\Omega)$ and consider
\begin{eqnarray*}
J &=& \ds\int_\Omega \nabla (\beta I - T_D)u_\e \cdot \nabla v
\\
&=&
\ds\int_{\Omega \setminus D} \beta \nabla u_\e \cdot \nabla v
\;+\; 
\ds\int_{D} (\beta -1) \nabla u_\e \cdot \nabla v,
\end{eqnarray*}
in view of the definition of $T_D$. 
Inserting the expression~(\ref{def_ue}) of $u_\e$, we see that
\begin{eqnarray*}
J &=&
s_\e \ds\int_{\Omega \setminus D} \beta \nabla u \cdot \nabla (\chi_1^\e \chi_2 v)
\;+\; 
s_\e\ds\int_{D} (\beta -1) \nabla u \cdot \nabla (\chi_1^\e \chi_2 v)
\\
&&
\;+\; 
s_\e \ds\int_{\Omega \setminus D} \beta u \nabla(\chi_1^\e \chi_2) \cdot \nabla v
\;+\;
s_\e\ds\int_{D} (\beta -1) u \nabla(\chi_1^\e\chi_2) \cdot \nabla v
\\
&&
\;-\;
s_\e \ds\int_{\Omega \setminus D} \beta \nabla u \cdot v \nabla(\chi_1^\e\chi_2)
\;-\;
s_\e\ds\int_{D} (\beta -1) \nabla u \cdot v \nabla(\chi_1^\e\chi_2).
\end{eqnarray*}
Since $u$ is a local solution to~(\ref{eq_R2}), the sum of the first 2 integrals
vanishes, and we remain with
\begin{eqnarray} \label{est_J}
J &=&
\left( 
s_\e \ds\int_{\Omega} a u \nabla(\chi_1^\e \chi_2) \cdot \nabla v
\;+\;
s_\e\ds\int_{\Omega \cap (B_{2r_0} \setminus B_{r_0})} 
a u \nabla(\chi_2) \cdot \nabla v \right)
\nonumber \\
&&
\;+\;
s_\e\ds\int_{\Omega \cap (B_{2\e} \setminus B_\e)} 
a v \nabla(\chi_1^{\e}) \cdot \nabla u
\;=:\; s_\e(J_3+J_4).
\end{eqnarray}
where $a = \beta$ in $\Omega \setminus \overline{D}$ and $a = \beta -1$
in $D$.
The Cauchy-Schwarz inequality allows us to estimate the 
first two terms on the right-hand side by
\begin{eqnarray} \label{est_J3}
|J_3| &\leq& C s^\e \, ||v||_{H^1}
\left\{
\ds\int_0^{2\e} \ds\int_0^{2\pi}
\left(
|u|^2 |\chi_1^\prime|^2/\e^2 + |\nabla u|^2 |\chi_1|^2 
\right)\,rdrd\theta
\right.
\nonumber \\
&& \;+\; 
\left.
\ds\int_{r_0}^{2r_0} \ds\int_0^{2\pi}
\left(
|u|^2 |\chi_2^\prime|^2 + |\nabla u|^2 |\chi_2|^2 
\right)\,rdrd\theta
\right\}.
\end{eqnarray}
and the same arguments as those used to control the term $J_1$ in~(\ref{est_J1})
show that the two integrals above are $O(1)$.
As for the last term in~(\ref{est_J}), we write 
\begin{eqnarray*}
J_4 &:=& \ds\int_{B_{2\e} \setminus B_\e}
a \nabla u \cdot v \nabla \chi_\e
\nonumber \\
&=&
\ds\int_{B_{2\e} \setminus B_\e} 
a \nabla u \cdot \overline{v} \nabla \chi_\e
\;+\;
\ds\int_{B_{2\e} \setminus B_\e} 
a \nabla u \cdot (v-\overline{v}) \nabla \chi_\e,
\end{eqnarray*}
where $\overline{v} = |B_{2\e}|^{-1} \int_{B_{2\e}} v(x) \,dx$.
We note that the first integral in the above right-hand side reduces to
\begin{eqnarray*}
\overline{v} \int_0^{2\pi}  a(\theta) \phi(\theta)\,d\theta
\int_{\e}^{2\e}
i \xi \, \left(\ds\frac{r^{i\xi -1}-r^{-i\xi-1}}{2} \right)
\ds\frac{\chi_1^\prime(r/\e)}{\e} \,r dr
&=& 0.
\end{eqnarray*} 
Indeed, since $\phi$ is a solution to 
$(a(\theta)\phi^{\prime}(\theta))^\prime - \xi^2 a(\theta) \phi(\theta) = 0$, 
with periodic boundary conditions,
it satisfies
\begin{eqnarray*}
\ds\int_0^{2\pi} a(\theta)  \phi(\theta) \, d \theta &=& 0.
\end{eqnarray*}
It follows that
\begin{eqnarray*}
|J_4| &\leq&
\left(
\ds\int_{B_{2\e} \setminus B_\e} a^2 |\nabla u \cdot \nabla \chi_\e|^2 \,dx
\right)^{1/2}
\left(
\ds\int_{B_{2\e}} |v - \overline{v}|^2 
\right)^{1/2}.
\end{eqnarray*}
Using the following Poincar\'e inequality
\begin{eqnarray*}
\ds\int_{B_{2\e}} |v - \overline{v}|^2
&\leq& 4 |B_{2\e}|^2 \ds\int_{B_{2\e}} |\nabla v|^2,
\end{eqnarray*}
we obtain
\begin{eqnarray*}
|J_4| &\leq&
C \e ||v||_{H^1(\Omega)} 
\left(
\ds\int_{0}^{2\pi} a(\theta)^2 |\phi(\theta)|^2 \,d\theta
\right)^{1/2}
\\
&& \quad\quad
\left( \ds\int_{\e}^{2\e} |i\xi \ds\frac{r^{i\xi -1}-r^{-i\xi -1}}{2}|^2 
\ds\frac{[\chi^\prime(r/\e)]^2}{\e^2}
\,r dr \right)^{1/2}
\\
&\leq&
C \left( \ds\int_{\e}^{2\e} r^{-1} dr \right)^{1/2}
||v||_{H^1(\Omega)}
\\
&\leq& C \sqrt{ln(2)}||v||_{H^1(\Omega)} 
\;=\; O(1) ||v||_{H^1(\Omega)}.
\end{eqnarray*}

Altogether, (\ref{est_J}, \ref{est_J3}) and the above estimate show that
\begin{eqnarray*}
\forall\; v \in H^1_0(\Omega),\quad
\left|
\ds\int_\Omega \nabla (\beta I - T_D)u_\e \cdot \nabla v
\right|
&\leq& O(s_\e) ||v||_{H^1},
\end{eqnarray*}
which proves the claim since $s_\e \to 0$.
\medskip

Finally, we show that $u_\e \to 0$ weakly in $H^1(\Omega)$. In fact, 
since this sequence is uniformly bounded in $H^1$, 
it suffices to show that $u_\e \to 0$ strongly in $L^2$, which
follows from~(\ref{def_ue}), from the boundedness of $\chi_1$ and $\chi_2$
and from the fact that $s_\e \to 0$.
\finproof

%%-----------------------------------------------------------------------

\section{Characterization of the essential spectrum}
%%-----------------------------------------------------------------------

In this section, we consider $\lambda \notin [\lambda_-, \lambda_+]$,
$\beta = 1/2 - \lambda$, and $k = (1 - 1/\beta)$. The latter satisfies
\begin{eqnarray}\label{hp_k}
k < k_+ = \frac{-(2\pi- \a)}{\a} < 0
&\textrm{or}&
k_- = \ds\frac{-\a}{2\pi-\a} < k < 0.
\end{eqnarray}
We show that $\beta \notin \sigma_{ess}(T_D)$, so that according to
Proposition~\ref{prop_spect}, $\lambda~\notin~\sigma_{ess}(\Ks)$.
\medskip

We proceed by contradiction: If $\beta \in \sigma_{ess}(T_D)$, then
there exists a singular Weyl sequence $u_\e$, that satisfies 
the conditions~(\ref{def_singW}).
In the next three sections, we show 
%%-----------------------------------------------------------------------
\begin{proposition} \label{thm_convWseq}
The sequence $u_\e$ converges to $0$ {\em strongly}
in $H^1(\Omega)$.
\end{proposition}
\medskip

%%-----------------------------------------------------------------------

This contradicts the fact that $||u_\e||_{H^1} = 1$.
Consequently, in view of Theorem~\ref{thm_singWseq}, this proves

%%-----------------------------------------------------------------------
\begin{theorem} \label{thm_mainresult}
The essential spectrum of $\Ks$ is exactly
\begin{eqnarray*}
\se(\Ks) &=& [\lm,\lp].
\end{eqnarray*}
\end{theorem}
\medskip

%%-----------------------------------------------------------------------

\subsection{Controling the energy of $u_\e$ away from the corner}
%%-----------------------------------------------------------------------

Let $z_\e = \beta u_\e - T_D u_\e \in H^1_0(\Omega)$.
Let $\rho < R_0$ and let $\chi_\rho$ denote a smooth, radial cut-off function,
such that
\begin{eqnarray*}
\chi_\rho(x) &=&
\left\{ \begin{array}{ll}
1 & \textrm{if}\; |x| \leq \rho/2,
\\
0 & \textrm{if}\; |x| \geq \rho.
\end{array} \right.
\end{eqnarray*}
Let $v_\e = (1-\chi_\rho) u_\e$. We show that

%%-----------------------------------------------------------------------
\begin{proposition} \label{prop_ve}
The sequence $v_\e$ converges strongly to $0$ in $H^1$.
\end{proposition}
%%-----------------------------------------------------------------------

{\bf Proof:} 
Assume that it is not the case. Then there exists $\d > 0$ and
a subsequence (still labeled with $\e$) such that
\begin{eqnarray} \label{est_ved}
||v_\e||_{H^1_0} &\geq& \d.
\end{eqnarray}
We note that for any $v \in H^1_0(\Omega)$,
\begin{eqnarray} \label{pde_ue}
\ds\int_\Omega \nabla z_\e \cdot \nabla v
&=&
\ds\int_\Omega \nabla (\beta u_\e - T_D u_\e) \cdot \nabla v
\\
&=&
\beta \ds\int_\Omega \nabla u_\e \cdot \nabla v
\;+\;
\ds\int_D \nabla u_\e \cdot \nabla v
\nonumber \\
&=&
\ds\int_\Omega a \nabla u_\e \cdot \nabla v,
\end{eqnarray}
where $a(x) = \beta, x \in \Omega \setminus \overline{D}$,
and $a(x) = \beta -1, x \in D$.
Given $v \in H^1_0(\Omega)$,  we compute
\begin{eqnarray*}
\ds\int_\Omega a \nabla v_\e \cdot \nabla v
&=&
\ds\int_\Omega a \nabla \left[ (1- \chi_\rho) u_\e \right] \cdot \nabla v
\\
&=&
\ds\int_\Omega a \left[ (1-\chi_\rho)\nabla u_\e - u_\e \nabla \chi_\rho \right]
\cdot \nabla v
\\
&=&
\ds\int_\Omega a \nabla u_\e \cdot
\left[ \nabla \left((1-\chi_\rho)v \right) + v \nabla \chi_\rho \right]
\;-\;
a u_\e \nabla \chi_\rho \cdot \nabla v
\\
&=&
\ds\int_\Omega
\nabla z_\e \cdot \nabla \left( (1-\chi_\rho)v \right)
\;-\; 
u_\e \nabla \cdot \left( a v \nabla \chi_\rho \right)
\;-\;
a u_\e \nabla \chi_\rho \cdot \nabla v.
\end{eqnarray*}
Invoking the Cauchy-Schwarz and the Poincar\'e inequality,
it follows that
\begin{eqnarray*}
\left|
\ds\int_\Omega \nabla \left((\beta I - T_D)v_\e\right) \cdot \nabla v
\right|
&=&
\left|
\ds\int_\Omega a \nabla v_\e \cdot \nabla v
\right|
\\
&\leq&
C\; \left( ||u_\e||_{L^2} + ||z_\e||_{H^1_0} \right)\, ||v||_{H^1_0},
\end{eqnarray*}
As $u_\e \to 0$ strongly in $L^2(\Omega)$ since it
converges weakly  to $0$ in $H^1$, we conclude that
\begin{eqnarray} \label{conv_Tve}
(\beta I - T_D)v_\e = (1-\chi_\rho) u_\e 
&\rightarrow& 
0 \quad \textrm{strongly in}\; H^1_0(\Omega).
\end{eqnarray}
We note that since $v_\e$ has support in $\Omega \setminus B_{\rho/2}$,
\begin{eqnarray*}
T_D v_\e &=& T_{\tilde{D}} v_\e,
\end{eqnarray*}
where $\tilde{D}$ denotes any {\em smooth} connected inclusion, such that
$(D \setminus B_{\rho/2}) \equiv (\tilde{D} \setminus B_{\rho/2})$,
and thus~(\ref{conv_Tve}) also reads 
\begin{eqnarray*} \label{conv_Ttildeve}
(\beta I - T_{\tilde{D}})v_\e = (1-\chi_\rho) u_\e 
&\rightarrow& 
0 \quad \textrm{strongly in}\; H^1_0(\Omega).
\end{eqnarray*}
It is easily seen that $v_\e \rightharpoonup 0$ weakly in $H^1_0(\Omega)$,
and, upon rescaling in view of~(\ref{est_ved}), we conclude from the
above estimate that $v_\e/||v_\e||_{H^1_0}$ is a singular Weyl sequence for $T_{\tilde{D}}$.
But $\Tilde{D}$ is smooth, so that the associated Neumann-Poincar\'e operator
is compact and does not have essential spectrum, which
contradicts this fact, and proves the Proposition.
\finproof

%%-----------------------------------------------------------------------

\subsection{Controling the energy of $u_\e$ near the corner}
%%-----------------------------------------------------------------------

We now focus on $w_\e := \chi_\rho u_\e$, which has compact
support in $B_\rho$.
In view of~(\ref{pde_ue}), 
it is easy to check that $w_\e$ satisfies
\begin{eqnarray*}
\partial^2_{rr} w_\e + 1/r \partial_r w_\e
+ 1/r^2 \partial^2_{\theta\theta} w_\e &=&
\tilde{f}_\e,
\end{eqnarray*}
in $\left( D \cap B_\rho \right)$ and in
$\left( (\Omega \setminus \overline{D}) \cap B_\rho \right)$.
The right-hand side is defined as
\begin{eqnarray*}
\tilde{f_\e} &=&
\chi_\rho \Delta z_\e 
+ b \nabla \chi_\rho \cdot \nabla u_\e
+ \nabla (b u_\e) \cdot \nabla \chi_\rho 
+ b u_\e \Delta \chi_\rho,
\end{eqnarray*}
and we note that it converges strongly to $0$ in $H^{-1}(\Omega)$.
Moreover, since the function $\chi_\rho$ is radial, 
$w_\e$ satisfies the following transmission conditions
on the edges of the corner
\[
\left\{ \begin{array}{lcl}
\;\; w_\e(r,\frac{\alpha}{2}|_{-}) \;=\; w_\e(r,\frac{\alpha}{2}|_{+}) 
\\
w_\e(r,-\frac{\alpha}{2}|_{-}) \;=\; w_\e(r,-\frac{\alpha}{2}|_{+}),
\\
\;\;(\beta-1) \partial_\theta w_\e(r,\frac{\alpha}{2}|_{-}) 
\;=\; \beta \partial_\theta w_\e(r,\frac{\alpha}{2}|_{+}), 
\\
(\beta-1) \partial_\theta w_\e(r,-\frac{\alpha}{2}|_{-}) 
\;=\; \beta \partial_\theta w_\e(r,-\frac{\alpha}{2}|_{+}),
\end{array} \right.
\]
where the notations $|_{-}$,  $|_{+}$ indicate taking the 
limit from left and right sides respectively. 

We set 
\begin{eqnarray} \label{def_A}
A &=& \ds\frac{\a}{2\pi-\a} \in (0,1),
\end{eqnarray} 
and consider the change of variables
$(r, \theta) \in (0,\rho) \times (-\a/2,\a/2)
\rightarrow (r, \pi - \theta/A)$, which maps $D \cap B_\rho$ into
$(\Omega \setminus \overline{D}) \cap B_\rho$.
We define
\[
\left\{ \begin{array}{lcl}
v_\e(r,\theta) &=& w_\e(r,\pi - \theta/A)
\\
\tilde{g}_\e(r,\theta) &=& \tilde{f}_\e(r, \pi - \theta/A),
\end{array} \right.
\quad\quad \textrm{for}\; (r,\theta) \in D \cap B_\rho.
\]
It is easy to check that when $(f,g) = (\tilde{f}_\e, \tilde{g}_\e)$,
the functions 
$(w,v) = (w_\e|_{D \cap B_\rho},v_\e)$ satisfy the following system
\begin{equation} \label{syst_weve}
\left\{ \begin{array}{lcl}
\partial^2_{rr} w + 1/r \partial_r w
+ 1/r^2 \partial^2_{\theta\theta} w &=& f,
\\
\partial^2_{rr} v + 1/r \partial_r v
+ A^2/r^2 \partial^2_{\theta\theta} v &=& g,
\end{array} \right.
\end{equation}
with the boundary conditions
\begin{equation} \label{bc_weve}
\left\{ \begin{array}{lcl}
v(\rho,\theta) &=& w(\rho,\theta) \;=\; 0,
\\
v(r,\pm \a/2) &=&  w(r,\pm \a/2),
\\
\partial_\theta v(r,\pm \a/2) &=&  
\ds\frac{-k}{A} \partial_\theta w(r,\pm \a/2).
\end{array} \right.
\end{equation}

In other words, $w_\e$ and $u_\e$ both satisfy an elliptic equation
and take nearly the same Cauchy data on the edges of the corner.
\medskip

To study the above system, we introduce the (closed) subspace 
$V_1 \subset H^1(D \cap B_\rho) \times H^1(D \cap B_\rho)$ of 
functions $(w,v)$ that satisfy
\[
\left\{ \begin{array}{lcll}
v(\rho,\theta) &=& w(\rho,\theta) \;=\; 0 & |\theta| < \a/2
\\
v(r,\pm \a/2) &=& w(r,\pm \a/2) & 0 < r < \rho.
\end{array} \right.
\]

%%-----------------------------------------------------------------------
\begin{theorem} \label{thm_elliptic}
The system~(\ref{syst_weve}--\ref{bc_weve}) has a unique solution
$(w,v) \in V_1$.
Moreover, there exists a constant $C > 0$, such that
\begin{eqnarray*} 
||\nabla w||_{L^2(D \cap B_\rho)} 
+ ||\nabla v||_{L^2(D \cap B_\rho)} 
&\leq& C\, \left( ||f||_{H^{-1}(D \cap B_\rho)} + 
||g||_{H^{-1}(D \cap B_\rho)} \right).
\end{eqnarray*}
\end{theorem}
%%-----------------------------------------------------------------------

{\bf Proof:}
On $V_1$ we consider the norm
\begin{eqnarray} \label{norm_V}
||(w,v)||
&:=&
\left(
\ds\int_{D \cap B_\rho} |\nabla w|^2 + |\nabla v|^2
\right)^{1/2}.
\end{eqnarray}
We multiply the equations~(\ref{syst_weve}) by two functions 
$\phi, \psi \in H^1(D \cap B_\rho)$ that vanish on $D \cap \partial B_\rho$, 
and integrate to obtain
\begin{eqnarray*}
\lefteqn{
\ds\int_{D \cap B_\rho} f \phi \;+\; g \psi}
\\
&=&
\ds\int_0^\rho \int_{-\a/2}^{\a/2}
\left( \begin{array}{c}
\partial_r w \\ \frac{1}{r} \partial_\theta w
\end{array} \right)
\cdot
\left( \begin{array}{c}
\partial_r \phi \\ \frac{1}{r} \partial_\theta \phi
\end{array} \right)
\;+\;
\left( \begin{array}{c}
\partial_r v \\ \frac{A^2}{r} \partial_\theta v
\end{array} \right)
\cdot
\left( \begin{array}{c}
\partial_r \psi \\ \frac{1}{r} \partial_\theta \psi
\end{array} \right)
\, r dr d\theta
\\
&&
\;-\;
\ds\int_{\theta = \pm \a/2}
\ds\frac{1}{r} \partial_\theta w \phi
\;+\;
\ds\frac{A^2}{r} \partial_\theta v \psi.
\end{eqnarray*}
We note that that the last integral can be rewritten as
\begin{eqnarray*}
\ds\int_{\theta = \pm \a/2}
\ds\frac{1}{r} \partial_\theta w (\phi - Ak \psi)
\;+\;
\ds\frac{A^2}{r} \left[
\partial_\theta v + \ds\frac{k}{A} \partial_\theta w
\right] \psi.
\end{eqnarray*}
To satisfy the natural boundary condition in~(\ref{bc_weve}),
we are thus led to introduce the subspace 
$V_2 \subset H^1(D \cap B_\rho) \times H^1(D \cap B_\rho)$ of 
functions $(\phi, \psi)$ that satisfy
\[
\left\{ \begin{array}{ll}
\phi(\rho,\theta) \;=\; \psi(\rho,\theta) \;=\; 0, & \quad |\theta| < \a/2
\\
\\
\phi(r,\pm \frac{\a}{2}) - Ak \psi(r,\pm \frac{\a}{2}) \;=\; 0, & \quad 0 < r < \rho,
\end{array} \right.
\]
which we also equip with the norm~(\ref{norm_V}).
We also introduce the following bilinear form $B$ on $V_1 \times V_2$ by
\begin{eqnarray*}
B\left(
\left( \begin{array}{c} w \\ v \end{array} \right),
\left( \begin{array}{c} \phi \\ \psi \end{array} \right)
\right)
&=&
\ds\int_0^\rho \int_{-\a/2}^{\a/2}
\left( \begin{array}{c}
\partial_r w \\ \frac{1}{r} \partial_\theta w
\end{array} \right)
\cdot
\left( \begin{array}{c}
\partial_r \phi \\ \frac{1}{r} \partial_\theta \phi
\end{array} \right)
\\
&&
\;+\;
\left( \begin{array}{c}
\partial_r v \\ \frac{A^2}{r} \partial_\theta v
\end{array} \right)
\cdot
\left( \begin{array}{c}
\partial_r \psi \\ \frac{1}{r} \partial_\theta \psi
\end{array} \right)
\, r dr d\theta.
\end{eqnarray*}
Thus, solving~(\ref{syst_weve}--\ref{bc_weve}) amounts
to solving the variational problem~:
find $W = (w,v) \in V_1$ such that
\begin{eqnarray*}
\forall\; \Phi = (\phi,\psi) \in V_2, \quad
B(W,\Phi)
&=& \ds\int_{D \cap B_\rho} f \phi + g \psi.
\end{eqnarray*}
It is easily checked that the above right-hand side
defines a continuous linear form on $V_2$, and that
\begin{eqnarray*}
\forall\; (W,\Phi) \in V_1 \times V_2,\quad
|B(W,\Phi)|
&\leq& ||W||\, ||\Phi||.
\end{eqnarray*}
Therefore, the theorem will be proved upon showing 
that $B$ satisfies the inf-sup condition~(for instance the version in~\cite{bab}),
i.e., that there exists $\delta > 0$ such that
\begin{eqnarray} \label{inf-sup}
\inf_{W \in V_1, ||W|| = 1} \left(
\sup_{\Phi \in V_2, ||\Phi||=1} 
B(W,\Phi) \right)
&\geq& \delta.
\end{eqnarray}
\medskip

%%-----------------------------------------------------------------------

Let $W = (w,v) \in V^1$ and $p, q, d \in \R$. We set
\begin{eqnarray*}
\phi \;=\; (Akp+d)w + (Akq-d)v, &\quad&
\psi \;=\; pw + qv,
\end{eqnarray*}

so that $\phi - Ak\psi  = d(w - v) = 0$, and thus
$(\phi,\psi) \in V_2$ is an admissible test function.
The integrand in the expression of $B(W,\phi)$ takes the form
\begin{eqnarray*}
e &:=&
\partial_r u 
\partial_r \left[(Akp+d)w + (Akq-d)v\right]
\;+\; 
\partial_r v \partial_r(pw+qv)
\\
&& \;+\;
r^{-2} \partial_\theta u
\partial_\theta \left[(Akp+d)w + (Akq-d)v\right]
\;+\; \ds\frac{A^2}{r^2} \partial_\theta v
\partial_\theta(pw+qv)
\\
&=&
(Akp+d) \xi_1^2 + (Akq-d+p)\xi_1\xi_3 + q \xi_3^2
\\
&&  
\;+\; 
(Akp+d) \xi_2^2 + (Akq-d+A^2p)\xi_2\xi_4 + A^2q \xi_4^2,
\end{eqnarray*}
where $\xi_1 = \partial_r \phi, \xi_2 = r^{-1} \partial_\theta \phi,
\xi_3 = \partial_r \psi, \xi_4 = r^{-1} \partial_\theta \psi$.
Fixing $q=1$, it follows that $e$ defines a positive definite quadratic form
(pointwise) provided that the polynomials
\begin{eqnarray*}
P_1(\xi) &=& (Akp + d) + (Ak - d + p)\xi + \xi^2,
\\
P_2(\xi) &=& (Akp + d) + (Ak - d + A^2p)\xi + A^2\xi^2,
\end{eqnarray*}
are strictly positive, in other words, provided that
\begin{equation} \label{discr}
\left\{ \begin{array}{lcl}
(Ak - d + p)^2 - 4(Akp + d) &<& 0,
\\
(Ak - d + A^2 p)^2 - 4A^2(Akp + d) &<& 0.
\end{array} \right.
\end{equation}
 We regard these expressions as polynomials in $p$, the roots of
 which are respectively
 \begin{eqnarray*}
 f_\pm(d ) &=& (d + Ak) \pm 2 \sqrt{d(1 + Ak)},
 \\
 g_\pm(d) &=& \ds\frac{1}{A^2}
 \left[ (d + Ak) \pm 2 A\sqrt{d(1 + k/A)} \right].
 \end{eqnarray*}
 We remark that the roots are real if and only if
 \[
 \left\{ \begin{array}{ll}
 k \;<\; -1/A &\quad \textrm{if}\; \lambda < 0,
 \\
 -A < k < 0 & \quad \textrm{if}\; \lambda > 0,
 \end{array} \right.
 \]
 i.e., recalling~(\ref{def_A}, \ref{hp_k}),
 if and only if $\lambda \notin [\lambda_-, \lambda_+]$,
 which is our hypothesis.
 \bigskip
 
It only remains to show that we can indeed find parameters 
$p,d$ for which~(\ref{discr}) is satisfied,
i.e. that we can find $d$ such that
\begin{eqnarray} \label{cond_intertw}
(f_-(d),f_+(d)) \cap (g_-(d), g_+(d))
&\neq& \emptyset
\end{eqnarray}
(and then pick $p$ in the intersection).
\medskip
 
To this end, assume first that $-A < k < 0$,
so that $d_+ := -Ak > 0$. We note that
\begin{eqnarray*}
\ds\frac{f_+(d_+) + f_-(d_+)}{2} &=& 0,
\end{eqnarray*}
and that
\begin{eqnarray*}
g_+(d_+) &=& \ds\frac{2}{A^2} \sqrt{Ak(k/A-1)} \;>0,
\\
g_-(d_+) &=& \ds\frac{-2}{A^2} \sqrt{Ak(k/A-1)} \;<0,
\end{eqnarray*}
which yields (\ref{cond_intertw}).
\medskip
 
If $k < -1/A$, one can see that $d_- = A^2 + kA <0$ and that
\begin{eqnarray*}
\ds\frac{f_+(d_-) + f_-(d_-)}{2} &=& 2Ak + A^2 \;<\; -1
\;=\; g_p(d_-).
\end{eqnarray*}
On the other hand, since $0 < A < 1$ and $k < -1/A$, we have
\begin{eqnarray*}
\ds\frac{-2}{A}\sqrt{d(1 + k/A)} &<& -2\sqrt{d(1+Ak)},
\end{eqnarray*}
so that for any $d < 0$, $g_-(d) < f_-(d)$, and in particular
$g_-(d_-) < \frac{f_+(d_-) + f_-(d_-)}{2}$.
It follows that~(\ref{cond_intertw}) also holds in this case.
\finproof

%%-----------------------------------------------------------------------

\subsection{Proof of Proposition~\ref{thm_convWseq}}
%%-----------------------------------------------------------------------

We come back to the singular Weyl sequence $u_\e$, which
we split as $u_\e = (1-\chi_\rho)u_\e + \chi_\rho u_\e$.
Proposition~\ref{prop_ve} shows that $(1-\chi_\rho)u_\e$ 
converges strongly to $0$.
On the other hand, Theorem~\label{elliptic} applied to
$\chi_\rho u_\e$ shows that
\begin{eqnarray*}
||\nabla (\chi_\rho u_\e)||_{L^2(B_\rho)}
&\leq&
C \, \left( 
||\tilde{f}_\e||_{H^{-1}(D \cap B_\rho)}
\;+\; 
||\tilde{g}_\e||_{H^{-1}(D \cap B_\rho)}
\right)
\\
&\leq&
C \, \left( ||z_\e||_{H^1(\Omega)} + ||u_\e||_{L^2(\Omega)} \right)
\;\rightarrow\;0.
\end{eqnarray*}
It thus follows that $u_\e$ converges strongly to $0$ in $H^1(\Omega)$, 
which contradicts the assumption that $||u_\e||_{H^1(\Omega)} = 1$,
so that $\beta \notin \sigma_{ess}(T_D)$.
\finproof
\medskip

%%-----------------------------------------------------------------------

\noindent \textbf{Acknowledgements} 
\\
The work of Hai Zhang was supported by HK RGC grant ECS 26301016 and 
startup fund R9355 from HKUST.
E. Bonnetier was partially supported by the AGIR-HOMONIM grant from 
Universit\'e Grenoble-Alpes, and by the Labex PERSYVAL-Lab (ANR-11-LABX-0025-01).
This project was initiated while E.B. was visiting Hong Kong University of
Science and Technology, and completed at the Institute
of Mathematics and its Applications at the University of Minnesota.
The hospitality and support of both institutions is gratefully acknowledged.

%%-----------------------------------------------------------------------

%%-----------------------------------------------------------------------
\end{document}